\documentclass[11pt]{amsart}
\usepackage{palatino}
\usepackage{graphicx}
\usepackage{amsmath,amsthm,amssymb}

\theoremstyle{definition}

\newtheorem*{definition*}{Definition}
\newtheorem*{example*}{Example}
\newtheorem*{remark*}{Remark}
\newtheorem*{question*}{Question}
\newtheorem*{problem*}{Problem}
\newtheorem*{note*}{Note}
\newtheorem*{claim*}{Claim}

\numberwithin{equation}{section} 

\newcommand{\R}{\mathbf{R}}

\newcommand{\pr}{\mathrm{P}}

\newcommand{\eps}{\varepsilon}

\title{From Bocce to Positivity: \\Some Probabilistic Linear Algebra}
\author{Kent E. Morrison}
\email{morrison@aimath.org}
\date{March 1, 2013 \\ \indent Published in \emph{Mathematics Magazine} 86 (2013) 110--119}
\address{American Institute of Mathematics\\360 Portage Ave.\\ Palo Alto, CA 94306}
\address{California Polytechnic State University\\San Luis Obispo, CA 93407}
\begin{document}

\begin{abstract}
A question in geometric probability about the location of the balls in a game of bocce leads to related questions about the probability that a system of linear equations has a positive solution and the probability that a random zero-sum game favors the row player. Under reasonable assumptions we are able to find these probabilities.
\end{abstract}
\maketitle 
\markright{\textsc{FROM BOCCE TO POSITIVITY}}
\large
\renewcommand{\baselinestretch}{1.2}   
\normalsize

\section*{The Questions}  
At the American Institute of Mathematics we often play a lunch time game of bocce in a nearby park. Each round begins with the winner of the previous round tossing the pallino, a small white marker ball, out onto the lawn somewhere. Then each player takes turns throwing larger balls as close as possible to the pallino. Winners like to think their good fortune is a result of skill, but we see a lot of variability in the results leading us to conclude that chance plays a major role. Usually the balls cluster around the pallino, but from time to time they do not, and once, when all eight of them were off to one side---so that the pallino was not contained in their convex hull---I became intrigued with the possibility of figuring out the probability of that occurring.

Answering this question leads naturally into other questions that can be answered with the same approach. One of these questions is about positive solutions of systems of linear equations. We might be especially interested in positive solutions because the equations are a model for some real world problem where negative values of the variables do not make sense. Suppose that a system of liner equations is drawn out of a hat---that is, the coefficients of the system are random. What is the probability that there is a positive solution?

A third question concerns random two-person zero-sum games. The payoff matrix for such a game is an $m \times n$ rectangular matrix of real numbers. One player chooses a row, the other player a column, and the column player pays the row player the amount of the corresponding matrix entry---a negative number meaning that the row player actually pays the column player. Most of the time players do not have unique best choices, but they do have optimal probabilistic strategies specifying the probabilities for choosing among their options. 
It stands to reason that the player having more options is more likely to have the advantage, but can we quantify that? More precisely, assuming the game matrix has random entries as likely to be positive as negative, what is the probability that the row player expects a positive payoff?
\section*{Some Low Dimensional Situations}

For each question some special cases are easy to answer, and we begin with those. 

(1) Bocce as it is played in this world is a two-dimensional game, but imagine playing bocce in $d$ dimensions for any $d \geq 1$. The special case of $d=1$ we'll take care of now.  Place the pallino at the origin on the real line, and assume that the $n$ players' balls are independent random points with equal probability to be on either side of the origin and zero probability to be exactly at the origin. Then the origin fails to be in the convex hull of the random points exactly when all the points are on the same side of the origin. The first point can be on either side with the remaining $n-1$ points on the same side as the first point, an event that has probability $1/2^{n-1}$.
\vspace{3mm}

(2) Let's start with a single linear equation in two variables \[a_1 x_1 + a_2 x_2 = b,\] so that the set of solutions is a line in the plane---as long as $a_1$ and $a_2$ are not both 0. Assume that the $a_1$ and $a_2$ are independent random variables each having probability $1/2$ of being positive, probability $1/2$ of being negative, and probability zero of being 0. When $b \neq 0$, there are two intercepts, $x_1 = b/a_1$ and $x_2 = b/a_2$, and four possibilities for the signs of the intercepts, each having probability $1/4$. Unless both intercepts are negative, the line will meet the first quadrant and there will be positive solutions. Therefore the probability of a positive solution is $3/4$. However, when $b=0$ the solution set is a line through the origin with slope $-a_1/a_2$. The sign of the slope is equally likely to be positive as negative, and so the probability is $1/2$ that the line meets the first quadrant, giving a positive solution.

A similar analysis works for one equation in $n$ variables
\[ a_1 x_1 + \cdots + a_n x_n = b. \]
For $b \neq 0$ the probability of a positive solution is $1 - 1/2^{n}$, while for $b = 0$ the probability is $1 - 1/2^{n-1}$.

(3) Let's consider the game with $m=1$ and arbitrary $n$. The row player has only one choice, and the column player just picks the smallest entry since he wants to minimize what he pays. The minimum value will be positive when all $n$ entries are positive. Hence, if we assume that the entries are independent and just as likely to be positive as negative and with no chance of being zero, then the probability that they are all positive is $1/2^n$, and that is the probability that the row player has the advantage. 

\section*{Two-Dimensional Bocce}
We return to the original bocce game as it is played on a two-dimensional surface. The pallino is at the origin and the players' balls are random points $z_1,\ldots,z_n$ in $\R^2$. In actual play, especially with skilled players, the locations of the balls are not independent, because players attempt to knock opponents' balls out of the way or to set up blocking positions in front of the pallino. However, our casual game is not played on a smooth bocce court but rather in a park with bumpy terrain, patches of dead lawn, trees, hills, sidewalks, and playground equipment. The player with the nearest ball on each round begins the next round by throwing the pallino wherever he or she chooses, and so all the irregularities of the terrain make it reasonable to assume that the location of the players' balls are independent random points. The second key assumption is that the probability distribution for each ball has a density function symmetric with respect to the origin. For a probability density function $f$ this means that $f(-z)=f(z)$ for $z \in \R^2$. (We are not requiring that the random points be identically distributed.) Common examples include bivariate normal distributions centered at the origin and uniform distributions on regions such as disks or rectangles centered at the origin. It follows from this assumption that for any line through the origin, the probability is $1/2$ to be in each of the open half-planes on either side of the line. 

With our assumptions in place let's define $E$ to be the event that the origin does not lie within the convex hull of $n$ random points $z_1,\ldots,z_n$. For such a configuration there is a unique distinguished point among the $n$ points with the property that all the remaining points are in the half-plane described by starting at that point and going $\pi$ radians counter-clockwise. This shows that $E$ is the disjoint union of events $E_1,\ldots,E_n$, where $E_i$ is the event that $z_i$ is the distinguished point. Therefore, 
\[\pr(E) = \sum_{i=1}^n \pr(E_i),\]
but $\pr(E_i) = 1/2^{n-1}$ because $z_i$ can be anywhere and the other $n-1$ points must be in the correct half-plane. Therefore,  \[ \pr(E) = \frac{n}{2^{n-1}} .\]

This result is consistent with the informal observations in our bocce games in the park. We usually have four players, each throwing two balls, and so $n=8$, in which case $\pr(E)=8/2^7=1/16$. 

\section*{Higher Dimensional Bocce}

The two-dimensional analysis does not seem to work in three dimensions because there is not a natural way to associate to each nonzero point in $\R^3$ a half-space with that point on the separating plane. The same goes for higher dimensions. But in 1962 J. G. Wendel \cite{Wendel62} found an elegant solution relying fundamentally on an old theorem of L. Schl\"afli that counts the number of regions in $\R^d$ created by $n$ generic hyperplanes through the origin. Schl\"afli's result is in \emph{Theorie der vielfachen Kontinuit\"at}, written between 1850 and 1852, which is one of the seminal contributions to the development of higher dimensional geometry in the nineteenth century. Despite repeated efforts by Schl\"afli and others it was almost 50 years before it was eventually published in 1901, six years after his death. Schl\"afli's result can be found in the free Google Books edition \cite[p. 41]{Schlafli01} or in his collected works \cite[p. 211]{Schlafli50}.

As Wendel puts it, there are $n$ points ``scattered at random on the surface of the unit sphere'' in $\R^d$, and the problem is to evaluate the probability that all the points lie on some hemisphere. But this is just what we want, because the origin is not in the convex hull of  nonzero points $z_1,\ldots,z_n$ in $\R^d$ if and only if the points all lie in some half-space, or equivalently that $z_1/|z_1|,\ldots,z_n/|z_n|$ all lie on some hemisphere of the unit sphere.

What follows is a modified version of Wendel's solution. Let $p(n,d)$ be the probability that the convex hull of $n$ random points in $\R^d$ does not contain the origin (equivalently, that the points lie in some half-space). We assume that the points $z_1,\ldots,z_n$ are independent and the probability distributions of the points are symmetric with respect to the origin and that they have density functions. For each of the $2^n$ sign vectors $\eps = (\eps_1,\ldots,\eps_n) \in \{\pm 1\}^n$, we define a random variable $X_\eps$. The value of $X_\eps$ is $1$ if the points $\eps_1 z_1,\ldots,\eps_n z_n$ all lie in a half-space; otherwise the value is $0$. Then $p(n,d)=E(X_{\mathbf{1}})$, where $\mathbf{1}=(1,\ldots,1)$. Because the distributions are symmetric with respect to the origin, $E(X_\eps)$ is independent of $\eps$, and so $2^n E(X_\mathbf{1})=\sum_\eps E(X_\eps)= E\big( \sum_\eps X_\eps \big)$. Therefore,
\[ p(n,d) = \frac{1}{2^n} E\big( \sum_\eps X_\eps \big).
\]

The next step is showing that the sum $\sum_\eps X_\eps$ is constant almost surely and that the value of the constant is the number of connected regions in $\R^d$ created by the $n$ hyperplanes through the origin that are orthogonal to the $z_i$.  Each of the regions corresponds to a sign vector $\eps$ which describes that region as a particular intersection of half-spaces. That is, two points $v$ and $w$ are in the same region if and only if the inner products $\langle v, z_i \rangle$ and  $\langle w, z_i \rangle$ have the same sign for $i=1,\ldots,n$.
But not all sign vectors correspond to regions because the intersection of half-spaces described by a sign vector can be empty. Now for each $\eps$ that does come from a region, let $v$ be a point in the region. Then $\eps_1 z_1, \ldots,\eps_n z_n$ all lie in the half-space of points $x$ such that $\langle v, x \rangle > 0 $, and so $X_\eps= 1$ for the random points $z_1,\ldots,z_n$. Summing over $\eps$ we conclude that $\sum_\eps X_\eps$ is equal to the number of regions. (There are configurations of the points for which the sum $\sum_\eps X_\eps$ does not achieve this value but is something less. This occurs when there is some unexpected linear dependence among the $z_i$. To be precise, the exceptional configurations for $n \leq d$ are those for which $z_1,\ldots,z_n$ are linearly dependent; for $n > d$ they are those for which some $d$ of the points are linearly dependent. For example, if $z_1,z_2,z_3$ in $\R^3$ all lie in a plane, then the planes normal to them divide space into only 6 regions rather than 8.) 

The last ingredient we need is Schl\"afli's formula for the number of regions. Letting $r(n,d)$ denote the number of regions created by $n$ hyperplanes through the origin in $\R^d$, what Schl\"afli proved is that
\[ r(n,d) = 2 \sum_{j=0}^{d-1} {n-1 \choose j} .\]
How could you discover this formula? Well,  $r(n,d)$ satisfies the recurrence relation
\[ r(n,d) = r(n-1, d) + r(n-1, d-1) \]
 with boundary conditions
$ r(n,2) = 2n$ and $r(n,d) = 2^n$ for $n \leq d$. The boundary conditions are straightforward, but the recurrence relation is more subtle and is explained in the appendix. With this information you can compute several values of $r(n,d)$ and then hope to notice that the difference $r(n,d)-r(n,d-1)$ is always twice a binomial coefficient; in fact, it is $2{n-1 \choose d-1}$. From that you build the formula, which can then be proved rigorously by showing it satisfies the recurrence relation and boundary conditions. 

With this we have everything we need to see that 
\[ p(n,d)=\frac{1}{2^n}r(n,d)=\frac{1}{ 2^{n-1}}\sum_{j=0}^{d-1} {n-1 \choose j}.\]

\begin{table}[htdp]
\begin{center}\begin{tabular}{ccccccccccccc} \hline
1 & 2 & 3 & 4 & 5 & 6 & 7 & 8 & 9 & 10 & 11 & 12&13 \\
1.0&1.0&1.0& .875& .686& .500& .344& .227& .145 & .090& 
.055& .033&.019 \\ \hline\end{tabular} \vspace{1mm}\caption{Values of $p(n,3)$ for $1 \leq n \leq 13$}
\end{center}
\label{defaulttable}
\end{table}
\begin{figure}
\vspace{.2in}
\centerline {
\includegraphics[width=4in]{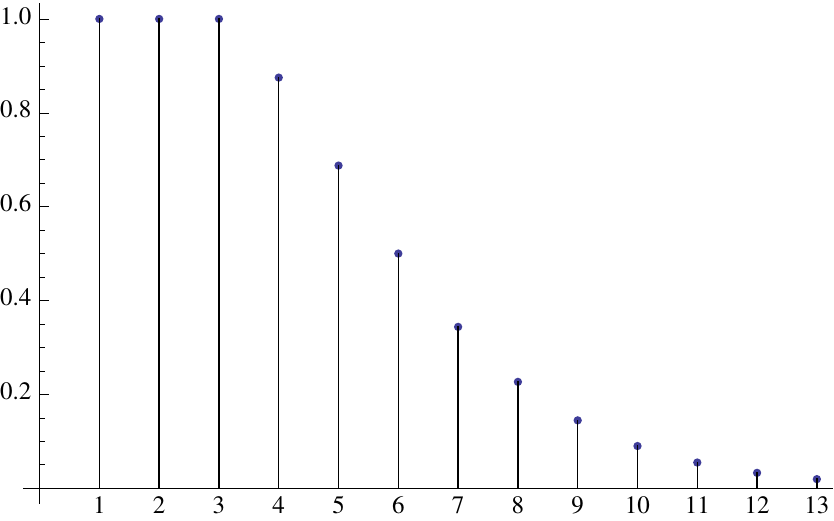}
}
\caption{The probability $p(n,3)$ for $1 \leq n \leq 13$.}
\vspace{.2in}
\end{figure}

With a little more work we can see that the location of the first ball does not really matter for the result. That is, you can specify the point $z_1$, while the other points are random, and the probability that they all lie in a half-space is still $p(n,d)$. To understand this, consider the sum of the $2^{n-1}$ random variables $X_\eps$ where $\eps$ ranges over the sign vectors in $\{\pm 1\}^n$ such that $\eps_1 = 1$. Then this sum is equal to the number of sign vectors corresponding to the regions that have $\eps_1 = 1$. But there is a one-to-one correspondence between those regions and the regions whose sign vectors have $\eps_1 = -1$, where the correspondence just pairs a region $R$ with its negative $-R$. Thus, the sum we want is $(1/2)r(n,d)$. Now we divide this by $2^{n-1}$ to get the expectation of $X_\mathbf{1}$, which is the probability that we want, but this gives us exactly $p(n,d)=(1/2^n)r(n,d)$.

 We end this section with the intriguing observation that 
\[ p(n,d)+p(n,n-d) = 1.\]
That is, $p(n,d)$ and $p(n,n-d)$ are complementary probabilities. 
This identity holds for $n \geq 0$ and all integer values of $d$ with the understanding that $p(n,d)=0$ for $d \leq 0$. Although it is easy to verify algebraically, we will return to it with a geometric proof in the section on random subspaces.
 
\section*{Positive Solutions of Linear Equations}
We can use what we already know about convex hulls in order to find the probability that a random system of linear equations has a positive solution. The key is that  
the existence of a positive solution for a homogeneous system of linear equations is equivalent to the property that the origin is in the convex hull of the column vectors of the  coefficient matrix. Here are the details. 

Consider an $m \times n$ matrix $A$ and the system of linear equations written as $Ax = 0$. If $x$ is a solution and $z_1, \ldots,z_n \in \R^m$ are the columns of $A$, then $\sum x_i z_i = 0 \in \R^m$.  Furthermore, if $x$ is a \emph{positive} vector (meaning that $x_i \geq 0$ for all $i$ and at least one of the $x_i$ is positive), then we can multiply $x$ by the scalar $1/\sum x_i$ to get a solution vector $t$ with the property that $t_i \geq 0$ and $\sum t_i = 1$, and thus $0$ is in the convex hull of the $z_i$. 

Therefore, if the columns of the matrix $A$ are random points in $\R^m$ that satisfy the assumptions in the bocce problem, then the probability that $Ax=0$ has a positive solution is the probability that $n$ points in $\R^m$ contain the origin in their convex hull, which is $1-p(n,m) = p(n,n-m)$. Those assumptions are satisfied if the entries of $A$ are independent random variables distributed with probability densities that are even functions. Examples include normal distributions and uniform distributions in balanced intervals of the form $[-c,c]$. 

Turning to the random system of equations $Ax=b$ where $b \neq 0$, we may consider $b$ as a random vector or as a fixed vector, because both lead to the same probability that there is a positive solution. Now suppose there is a positive solution $x$, and let $z_1,\ldots,z_n$ be the columns of $A$. Then $\sum x_i z_i = b$. Move $b$ to the other side and scale by $(1 + \sum x_i)^{-1}$ to see that $0$ is a convex combination of $-b$ and the $z_i$. Now we are in the situation of $n+1$ points with one fixed and the rest random, and in the previous section we determined that the probability that they lie in a half-space is $p(n+1,m)$. Therefore, the complementary probability $1-p(n+1,m) = p(n+1, n+1-m)$ is the probability that the origin is in the convex hull. When the origin is in the convex hull we have an expression of the form 
\[  -t_0 b + t_1 z_1 + \cdots + t_n z_n =0, \]
with $t_i \geq 0$, and we can move $b$ back to the right side to show a positive solution for $Ax = b$ as long as $t_0 \neq 0$. 

Let's consider the possibility that $t_0 =0$ in such a convex combination. Then $0$ is in the convex hull of the $z_i$ alone. Because any $m$ of the $z_i$ are linearly independent and $n > m$, the convex hull does not lie in a lower dimensional subspace and so it has non-empty interior. Since the probability is zero that the origin is on the boundary of the convex hull, it must be in the interior and so there is an open ball containing the origin and lying within the convex hull. For a sufficiently small $\lambda > 0$, the point $\lambda b$ is in the ball and hence in the convex hull of the $z_i$. Therefore, there exist $s_i \geq 0$ such that $ s_1 z_1 + \ldots + s_n z_n =\lambda b $. Multiplying both sides by $1/\lambda$ gives a positive solution of $Ax=b$. 

In summary, the probability of a positive solution of $Ax=0$ is $p(n,n-m)$. The probability of a positive solution of $Ax=b$, where $b$ is a fixed or random nonzero vector, is $p(n+1,n+1-m)$. As a partial check let's look back at the case of one equation and two unknowns ($m=1$ and $n=2$) for which $p(n,n-m)=p(2,1)=1/2$ and $p(n+1,n+1-m)=p(3,2)=3/4$. Those values are just what we found earlier. In general, just as in this special case, it is more likely that there is a positive solution when $b \neq 0$.

\section*{Random Subspaces}
Since the solutions of the equation $Ax=0$ form the null space of the matrix $A$, the question about positive solutions becomes a question about the probability that a random subspace of $\R^n$ contains a positive vector. We can assume that the subspaces have a fixed dimension $n-m$, because the assumptions about the distribution of the entries of $A$ imply that $\textrm{rank}(A) =m$, almost surely, so that the null space has dimension $n-m$.

The orthogonal complement to the null space of $A$ is the row space of $A$ or the range of $A^t$, and so it is a subspace of dimension $m$. That is, a random $m$-dimensional subspace of $\R^n$, for $m < n$, is generated by taking the span of $m$ random vectors in $\R^n$. With reasonable assumptions the vectors will be linearly independent with probability one, thus giving a span of dimension $m$.

What is the probability that a random $m$-dimensional subspace of $\R^n$ contains a positive vector? We'll see that the answer is $p(n,m)$. Let the random subspace be the row space of the random $m \times n$ matrix $A$ as before (the columns are independent random vectors having probability densities symmetric with respect to the origin). In particular, the individual entries in $A$ could be independent random variables whose probability density functions are even. 

Now we apply the theorem from linear algebra known as Gordan's Theorem of the Alternative \cite{Roman08}: a subspace $V$ of $\R^n$ contains a positive vector if and only if its orthogonal complement $V^\perp$ does not contain a strictly positive vector. 
Let $V$ be the row space of $A$ so that $V^\perp$ is the null space of $A$. Then the probability that $V^\perp$ does not contain a strictly positive vector is the same as the probability that it does not contain a positive vector, namely, $1-p(n,n-m)$, which is $p(n,m)$.

This, then, gives us a geometric explanation for the complementary probability identity mentioned earlier
\[ p(n,m)+p(n,n-m) = 1 .\]
We can understand the result as a \emph{probabilistic theorem of the alternative}: with probability one, a random subspace either contains a positive vector or its orthogonal complement contains a positive vector, but not both. In short, complementary subspaces define complementary events. 

\begin{figure}
\centerline {
\includegraphics[width=3in]{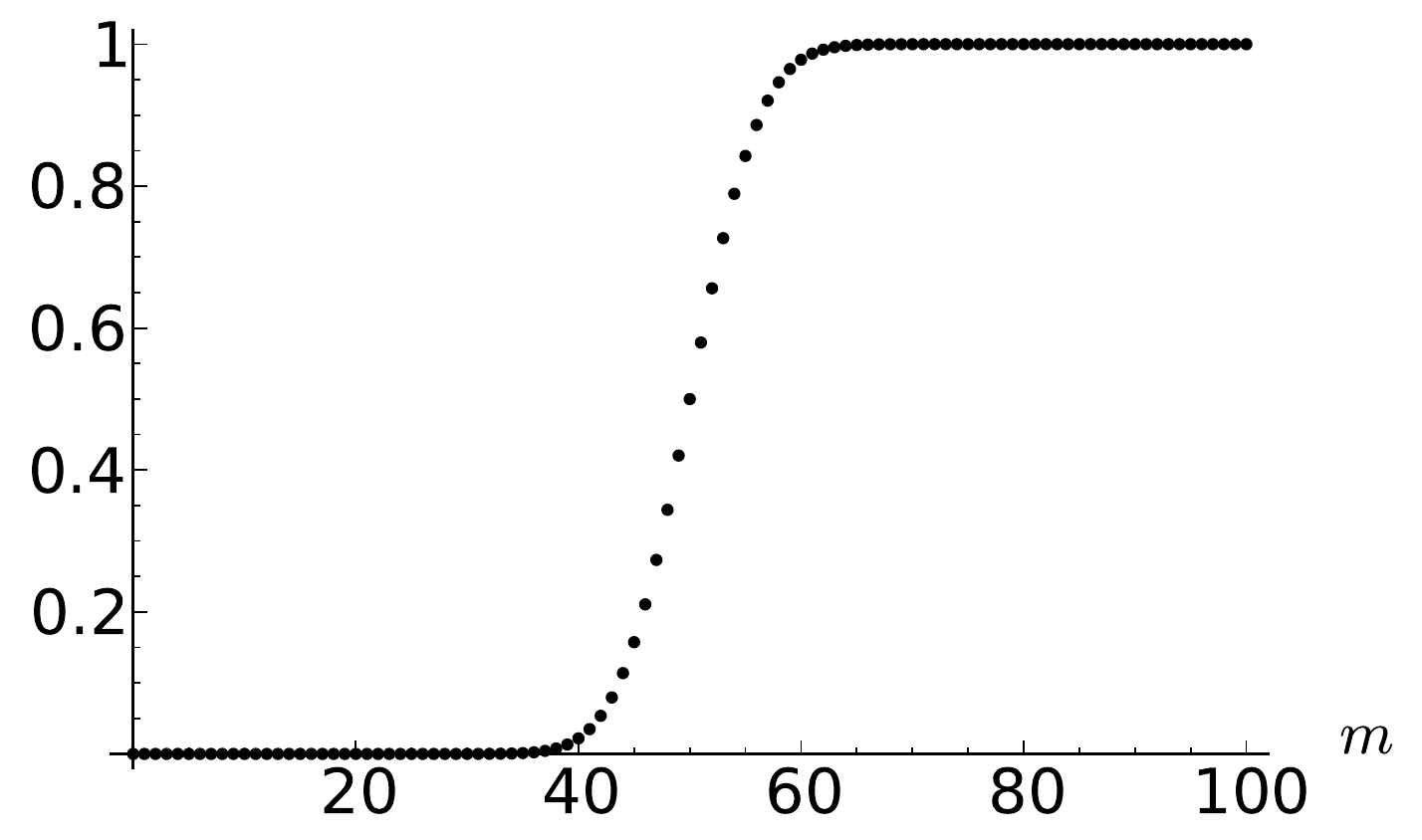} }
\caption{$p(n,m)$  for $n=100$.}
\end{figure}

Figure 2 shows that for large $n$ the value of $p(n,m)$ rises quickly from nearly 0 to nearly 1 as $m$ passes through $n/2$. As $n$ increases the plot  becomes more and more like a step function. Note that for $n$ even $p(n,n/2)$ is always $1/2$.

We end this section with an explanation for the individual binomial probabilities 
\[\frac{1}{2^{n-1}} { {n-1} \choose j } \] 
that we see in the sums for $p(n,m)$. This, of course, is the probability of exactly $j$ heads occurring in a sequence of $n-1$ tosses of a fair coin. Now instead of flipping a coin, we generate a sequence of $n$ independent random vectors $v_1,\ldots,v_n$ in $\R^n$. Let $V_m$ be the span of the first $m$ vectors. Then $V_m$ is a random subspace of dimension $m$. The probability that $V_m$ contains a positive vector while $V_{m-1}$ does not is 
\[p(n,m)-p(n,m-1)=\frac{1}{2^{n-1}}{ {n-1} \choose {m-1}}. \]
In other words, the integer valued random variable whose value is the least $m$ such that $V_m$ contains a positive vector has the same distribution as the number of heads in $n-1$ tosses of a fair coin.

\section*{Random Games}
Consider a random two-person zero-sum game described by an $m \times n$ matrix $A$. The rows represent the pure strategies of the row player and the columns the pure strategies of the column player, and the convention is that when the row player chooses row $i$ and the column player chooses column $j$, the result is that the column player pays the row player the amount $a_{ij}$. (If $a_{ij}$ is negative, then the column player receives $-a_{ij}$ from the row player.) Thus, positive entries are good for the row player and negative entries are good for the column player. We assume that the entries are independent random variables with probability density functions that are even, so that each individual entry favors neither player. 

In general there is not an optimal pure strategy for each player, but there are optimal mixed strategies. A mixed strategy is a probability distribution on the finite set of pure strategies. For the row player it is a probability vector $p=(p_1,\ldots,p_m)$ where $p_i$ is the probability of choosing row $i$ to play, and for the column player a mixed strategy is a probability vector $q=(q_1,\ldots,q_n)$ where $q_j$ is the probability of playing column $j$. The players choose their strategies independently and so the expected payoff to the row player is  
\[  \sum_{i,j} a_{ij}p_i q_j ,\]
which can be written as the product 
$p^t A q$, where $p$ and $q$ are treated as column vectors. 

The row player wants to choose $p$ to make this product as large as possible, while the column player wants to choose $q$ to minimize it. The Minimax Theorem, von Neumann's fundamental result, asserts that 
\[ \max_p \min_q p^tAq = \min_q \max_p p^t A q, \]
and this number is called the \emph{value} of the game. Furthermore, the theorem asserts that there exist \emph{optimal} mixed strategies $p_*$ and $q_*$, not necessarily unique, such that ${p_*}^t A q_*$ equals the value of the game. (For full treatment of this material see the book by Guillermo Owen \cite{Owen95} or the e-book by Tom Ferguson \cite{Ferguson}.)

The game favors the row player when the value is positive, since the value of the game is the expected amount that the row player receives when the players use their optimal strategies.
We'd like to know the probability of that event for a random $m$ by $n$ game. Intuition suggest that the game is more likely to favor the player with the greater number of strategies, and that for $m=n$ it should be equally likely that the game value is positive or negative.

\begin{figure} 
\centerline {
\includegraphics[width=3in]{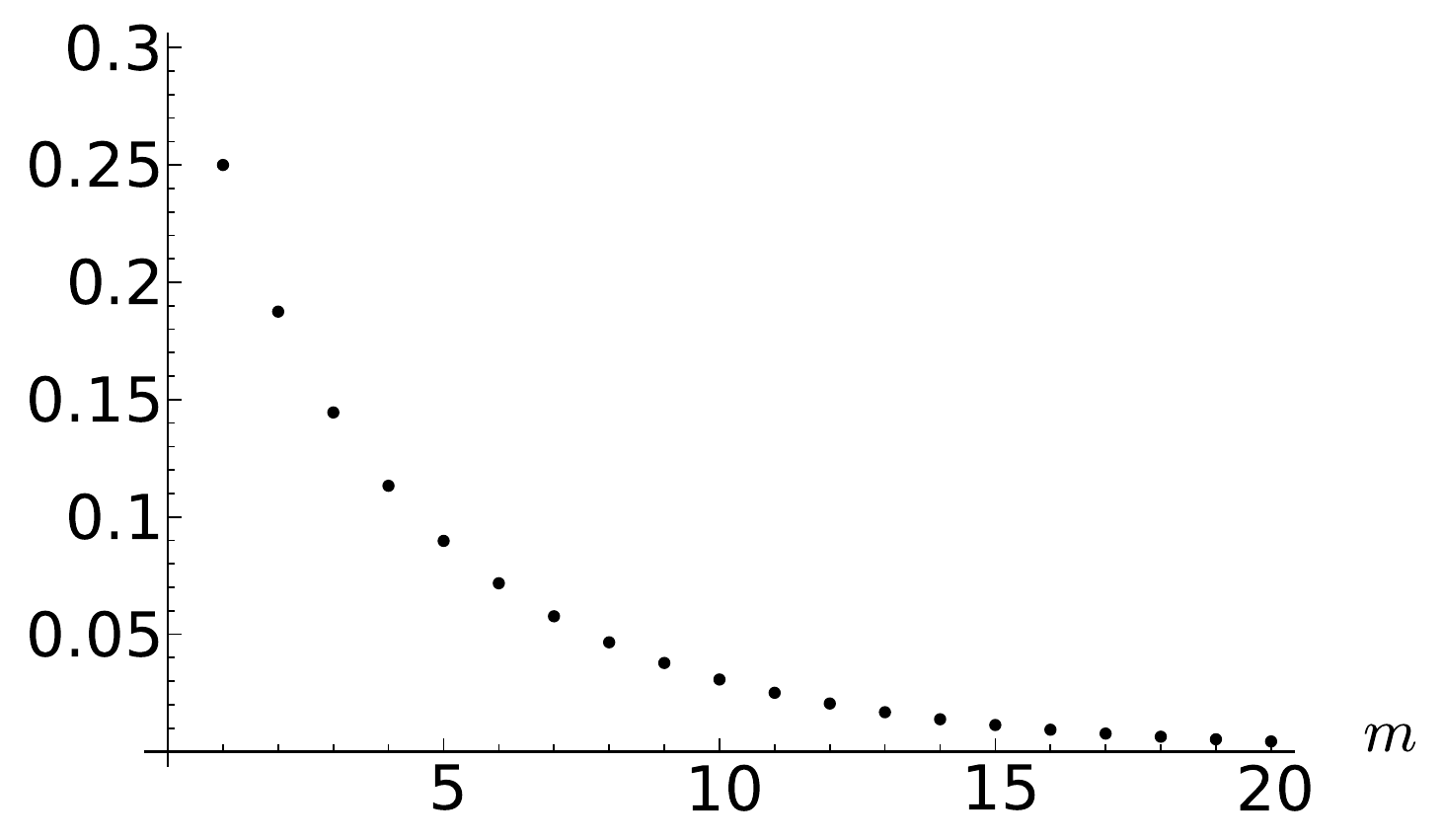}
}
\caption{Probability that row player wins with $m$ strategies when the column player has $2m$ strategies. \label{game-prob}}

\end{figure}

In 1966 Thomas Cover \cite{Cover66} proved---with reasonable assumptions on the entries of the payoff matrix---that the probability that the game favors the row player is $p(m+n,m)$. To get an idea of how much the advantage is for the player with more strategies, consider the situation in which the column player has twice as many pure strategies as the row player, i.e., $n=2m$. Figure \ref{game-prob} plots $p(n+m,m)$ vs. $m$, showing how quickly this probability approaches $0$ as $m$ increases.

To prove this result we assume that the entries of the payoff matrix $A$ are independent random variables whose probability density functions are symmetric with respect to 0.
In order for the value of the game to be positive, the row player must have a mixed strategy that gives a positive payoff no matter which column the opponent chooses. That is equivalent to the existence of a positive vector $x \in \R^m$ such that $\langle x,z_j \rangle >0$ for $j=1,\ldots,n$ where $z_j \in \R^m$ is the $j$-th column of the payoff matrix $A$. (From such an $x$ we get the mixed strategy by scaling it to make it a probability vector.)

If we multiply any of the rows and columns of the random $A$ by $-1$ the probability that the resulting matrix favors the row player does not change because of the symmetry of the matrix entries. Changing the signs of rows is the same as multiplying $A$ on the left by a diagonal matrix $C$ with diagonal entries $\pm 1$, while changing the signs of columns is equivalent to multiplying on the right by a diagonal matrix $D$ with $\pm 1$ entries. 

Let $W_{C,D}$ be the random variable equal to $1$ if the game with payoff matrix $CAD$ has positive value (i.e., favors the row player) and $0$ otherwise. The probability we want is then the expectation $E(W_{I,I})$, and  we have just noted that $E(W_{C,D})$ is independent of $C$ and $D$. 

Next we show that the sum of $W_{C,D}$ is constant (with probability one) where $C$ and $D$ range over all pairs of diagonal $\pm 1$ matrices. Consider the $m$ coordinate hyperplanes (i.e., the planes orthogonal to the standard basis vectors) along with the $n$ hyperplanes orthogonal to the $z_j$. These $m+n$ hyperplanes separate $\R^m$ into 
$r(m+n,m)$ regions by Schl\"afli's result. For $x$ in one of the regions let $c_i$ be the sign ($\pm 1$) of $x_i$ and let $d_j$ be the sign of $\langle x,z_j \rangle$. Let $C$ and $D$ be the diagonal matrices with diagonal entries $c_i$ and $d_j$. Then $x^t C$ is a strictly positive vector whose inner product with each column of $CAD$ is positive, and thus the payoff matrix $CAD$ describes a game that favors the row player.  
Therefore, with probability one, the random variable $\sum W_{C,D}$ is equal to the constant $r(m+n,m)$, and so $E(\sum W_{C,D}) =r(m+n,m)$. But
$ E(\sum W_{C,D}) = \sum E(W_{C,D}) = 2^{m+n} E(W_{I,I})$, and from this it follows that
\[ E(W_{I,I})=\frac{1}{2^{m+n}}r(m+n,m)= p(m+n,m). \]

\section*{Positive Input, Positive Output}
A linear system described by the $m \times n$ matrix $A$ maps an input vector $x$ in $\R^n$ to an output vector $Ax$ in $\R^m$. Again the linear system may only make sense as a physical model when there exists some positive input $x$ whose resulting output $Ax$ is also positive. Now suppose that $A$ is a random matrix with the same conditions on the independent entries as before. What is the probability that there is a positive input with positive output?

It turns out that we already have the answer to that question from the game theory situation, although in transposed form. Asking for a strictly positive output, rather than just a positive ouput, doesn't change the probability. For $Ax$ to be strictly positive means that the inner product of $x$ with each row of $A$ is positive. If $x$ is also positive, then we have the conditions for the game with matrix $A^t$ to have a positive value, and so we conclude that the probability of that occurring is $p(n+m,n)$. 

An easy case to check directly is $m=1$. Then the linear system is given by the map $x \mapsto \sum a_i x_i$. As long as any $a_i > 0$, it is possible to make $x_i$ a large positive number and the other $x_j$ small positive numbers so that the sum is positive. The complementary event that all the $a_i$ are negative has probability $1/2^n$, and so the probability that there is a positive input with positive output is $1-1/2^n$, which is indeed equal to $p(n+1,n)$.

\section*{The Recurrence Relation for $r(n,d)$}
This follows Wendel's \cite{Wendel62} paraphrase of Schl\"afli's proof.
The recurrence formula comes from analyzing how the number of regions changes as a hyperplane is added to the system. 
Let $H_1,\ldots,H_n$ be hyperplanes through the origin in $\R^d$. These hyperplanes are in general position meaning that any intersection of $k \leq d$ of them is a subspace of dimension $d-k$.  Omit $H_n$ for a moment and consider the regions created by the remaining $n-1$ hyperplanes. There are $r(n-1,d)$ regions, which are of two types---those that meet $H_n$ and those that don't. Let $\tau_1$ and $\tau_2$ be the number of each type; thus $r(n-1,d)=\tau_1 + \tau_2$. Now restore $H_n$ to the system. It cuts each region of type 1 into two parts, and so $r(n,d)= 2 \tau_1 + \tau_2$. Therefore, $r(n,d)=r(n-1,d)+\tau_1$. 

Now $\tau_1$ is also the number of regions in $H_n \cong \R^{d-1}$ created by the $n-1$ hyperplanes $H_i \cap H_n$, $i=1,\dots,n-1$, and so $\tau_1 = r(n-1,d-1)$. Therefore,
$ r(n,d) = r(n-1,d) + r(n-1,d-1).$

\end{document}